\documentclass[a4paper,11pt]{article}
\usepackage[hidelinks]{hyperref}
\hypersetup{
    colorlinks=true,
    linktoc=all,
    linkcolor=red,     
    citecolor=blue,     
}
\usepackage{inputenc}
\usepackage{graphicx}
\usepackage{amssymb}
\usepackage{latexsym, bm}
\usepackage{multicol}
\usepackage{indentfirst}
\usepackage{amssymb,amsfonts}
\usepackage{amsmath}
\usepackage{setspace}
\usepackage{enumerate}
\newtheorem{theorem}{Theorem}[section]
\newtheorem{lemma}{Lemma}[section]

\newtheorem{conjecture}{Conjecture}[section]

\usepackage{amssymb}
\newcommand{\ignore}[1]{}

\date{}

\begin{document}

\begin{spacing}{1}

\title{The Ramsey number of the 4-cycle versus a book graph}
\author{{Chunyang Dou,\footnote{School of Mathematical Sciences, Anhui University,
Hefei  230601, P.~R.~China. Email:{\tt chunyang@stu.ahu.edu.cn}.}}
~~~
{Tianyu Li,\footnote{College of Mathematics and Information Technology, Hebei Normal University Of Science \& Technology,
Qinhuangdao, 066004, P.~R.~China. Email: {\tt lty4765382@126.com.}}}
~~
{Qizhong Lin,\footnote{Center for Discrete Mathematics, Fuzhou University,
Fuzhou, 350108, P.~R.~China. Email: {\tt linqizhong@fzu.edu.cn}. Supported in part  by National Key R\&D Program of China (Grant No. 2023YFA1010202).}}
~~
 {Xing Peng\footnote{Center for Pure Mathematics, School of Mathematical Sciences, Anhui University, Hefei 230601, P.~R.~China. Email: {\tt x2peng@ahu.edu.cn}. Supported by the NSFC grant
(No.\ 12471319).}}
}

\maketitle

\begin{abstract}
 Given positive integers $n$ and $k$,  the  book graph $B_n^{(k)}$  consists of $n$ copies of $K_{k+1}$ sharing a common $K_k$.  The  book graph is a common generalization of a star and a clique,  which can be seen by taking $k=1$ and $n=1$ respectively.  In addition, the Ramsey number of a book graph is closely related to the diagonal Ramsey number. Thus the study of extremal problems related to  the book graph is of substantial  significance.
 In this paper,
 we aim to investigate the Ramsey number  $r(C_4,B_n^{(k)})$ which is the smallest integer $N$ such that for any graph $G$ on $N$ vertices, either $G$ contains $C_4$ as a subgraph or the complement $\overline{G}$ contains $B_n^{(k)}$ as a subgraph.  For $k=1$,   a  pioneer work by Parsons  ({\it Trans.~Amer.~Math.~Soc.,} 209 (1975), 33--44)  gives an upper bound for  $r(C_4,B_n^{(1)})$, which is tight for infinitely many $n$.  For $k=2$, in a recent paper ({\em J. Graph Theory,}  103 (2023), 309--322),  the second, the third, and the fourth authors  obtained the exact value of  $r(C_4,B_{n}^{(2)})$ for infinitely many $n$.
  The goal of this paper is to prove a similar result for each integer $k \geq 3$.  To be precise,  given  an integer $k \geq 3$ and a constant  $0<\varepsilon<1$,  let   $n=q^2-kq+t+\binom{k}{2}-k$ and $Q(k,\varepsilon)=(320k^4)^{k+1}/\varepsilon^{2k}$, where $1 \leq t \leq (1-\varepsilon)q$. We  first establish an upper bound for  $r(C_4,B_n^{(k)})$ provided $q \geq Q(k,\varepsilon)$.  Then
  we show the upper bound is tight for  $q \geq Q(k,\varepsilon)$ being a  prime power and  $1 \leq t \leq (1-\varepsilon)q$  under some assumptions. The proof leverages on a simple but novel refinement of a well-known inequality related to a $C_4$-free graph. Therefore, for each $k \geq 3$, we obtain the exact value of $r(C_4,B_n^{(k)})$ for infinitely many $n$.  Moreover, we prove general upper and lower bounds of  $r(C_4,B_n^{(k)})$ for $k \geq 3$.


\end{abstract}

\section{Introduction}
For two graphs $H_1$ and $H_2$, the {\it Ramsey number} $r(H_1,H_2)$ is the smallest integer $N$ such that for any graph $G$ on $N$ vertices, either $G$ contains $H_1$ as a subgraph or the complement $\overline{G}$ contains $H_2$ as a subgraph.
 For a graph $G$, if  neither $G$  contains $H_1$ as a subgraph nor $\overline{G}$ contains $H_2$ as a subgraph, then
 $G$  is called an $(H_1,H_2)$-{\it Ramsey graph}.
 Given positive integers $n$ and $k$,  let $B_n^{(k)}$ be the book graph which consists of $n$ copies of $K_{k+1}$ all sharing a common $K_{k}$.  Note that
$B_n^{(1)}$ is the star $K_{1,n}$ and $B_1^{(k)}$ is the complete graph $K_{k+1}$.   Therefore, the  book graph $B_n^{(k)}$ is a common generalization of the star and the complete graph.

 In this paper, we
focus on the Ramsey number of $C_4$ versus $B_n^{(k)}$ which enjoys a vivid interest. For the special case  where $k=1$, Parsons \cite{TP75} showed the following celebrated result, where the upper bounds were obtained  by using the double counting method (see Lemma \ref{C4-eg}) while the lower bounds were proved  by analyzing  structures of the $(C_4,B_n^{(1)})$-Ramsey graphs and the well-known Erd\H{o}s-R\'enyi orthogonal polarity graph \cite{Br66,ER66}.
\begin{theorem}[Parsons \cite{TP75}]\label{pars}
For any integer $n\geq 2$,
 $r(C_4, B_n^{(1)})\leq n+\lfloor\sqrt{n-1}\rfloor+2$.
Moreover, if $n=k^2+1$, then
 $r(C_4, B_n^{(1)})\leq n+\lfloor\sqrt{n-1}\rfloor+1$. Both upper bounds are tight.
In particular, for any prime power $q$, $r(C_4, B_{q^2}^{(1)})=q^2+q+1$ and $r(C_4,B_{q^2+1}^{(1)})=q^2+q+2$.
\end{theorem}

Subsequently,  Parsons  \cite{TP76} initiated the study of $r(C_4, B_{n}^{(1)})$ in which $n=q^2-t$ with $t>0$. The result from \cite{TP76} together with the one by Zhang, Chen, and Cheng \cite{ZC17} give the exact value of   $r(C_4,B^{(1)}_{q^2-t})$, where
 $1 \leq t \leq q+1$ and $t \neq q$ for $q$ being an even prime power as well as  $1 \leq t \leq 2\lceil \tfrac{q}{4}\rceil$ and $t \neq 2\lceil \tfrac{q}{4}\rceil-1$ for $q$ being an odd prime power.  In addition, Zhang, Chen, and Cheng \cite{ZC17F}  obtained the exact value of $r(C_4,B^{(1)}_{(q-1)^2+t})$ for $q \geq 4$ being an even  prime power and $t \in \{-2,0,1\}$ as well as the exact value of $r(C_4,B^{(1)}_{q(q-1)-t})$ for $q \geq 5$ being an odd  prime power and $t \in \{2,4,\ldots,2\lfloor \tfrac{q}{4}\rfloor\}$. The main work to prove results above is to construct the desired $(C_4, B_n^{(1)})$-Ramsey graph based on the Erd\H{o}s-R\'enyi orthogonal polarity graph. For small $n$, the exact value of $r(C_4, B_{n}^{(1)})$ is known for $2 \leq n \leq 50$, see \cite{TP75,TP76,ZBC}.

 Note that all known values of the Ramsey number $r(C_4, B_{n}^{(1)})$ for $n \geq 6$ satisfy either $r(C_4, B_{n}^{(1)})=n+\lfloor\sqrt{n-1}\rfloor+1$ or $r(C_4, B_{n}^{(1)})=n+\lfloor\sqrt{n-1}\rfloor+2$. It is an open question, see \cite{ZC17}, to show whether it is always true that $r(C_4, B_{n}^{(1)})=n+\lfloor\sqrt{n-1}\rfloor+\gamma$ for $\gamma \in \{1,2\}$.
If  one can answer this question in a positive way, then it will disprove the following conjecture due to Burr, Erd\H{o}s, Faudree, Rousseau, and Schelp \cite{befrs}.
\begin{conjecture}
For any constant $c>0$, we have $r(C_4, B_{n}^{(1)})<  n+\sqrt{n}-c$ for  infinitely many $n$.
\end{conjecture}

 Erd\H{o}s, as one of the authors, offered a \$100 prize for a proof or disproof of the  conjecture above.
 It is known \cite{befrs} that  $r(C_4, B_n^{(1)}) > n+\lfloor\sqrt{n}-6n^{11/40}\rfloor$ for all sufficiently large $n$.

Although the research on the $r(C_4, B_n^{(1)})$ is  extensive, for a long time, the value of $r(C_4,B_n^{(k)})$ is much less known for $k \geq 2$.
Previously, the exact value of $r(C_4,B_n^{(2)})$  is  determined only for $n \leq 14$, see \cite{GS71,FRS78,Tse1,Tse2}. Faudree, Rousseau and Sheehan \cite{FRS78} established the following  upper bound for $r(C_4,B_n^{(2)})$,
\begin{equation}\label{frs}
r(C_4,B_n^{(2)})\le g(g(n)), \;\;\text{where}\;\;g(n)=n+\lfloor\sqrt{n-1}\rfloor+2.
\end{equation}
As a special case in which $q$ is a prime power, it was proved in \cite{FRS78} that
\begin{align*}
q^2+q+2\leq r(C_4,B_{q^2-q+1}^{(2)})\leq q^2+q+4.
\end{align*}
Let $\mathcal{G}(q)$ be the set of $(C_4, B_{q^2-q+1}^{(2)})$-Ramsey graphs $G$ on $q^2+q+3$ vertices. If $\mathcal{G}(q)\neq\emptyset$, then we have $r(C_4,B_{q^2-q+1}^{(2)})= q^2+q+4$ for a prime power $q$. We know $\mathcal{G}(2) =\emptyset$ since $r(C_4,B_3^{(2)})=9$, see \cite{Tse1}. Faudree, Rousseau and Sheehan \cite{FRS78} proved that $\mathcal{G}(3) \not =\emptyset$ and asked to determine whether $\mathcal{G}(q)\neq\emptyset$ for each prime power $q>3$.
The answer to this problem is negative for $q=4$ as $r(C_4,B_{13}^{(2)})=22$, see \cite{Tse2}, which implies  $\mathcal{G}(4)=\emptyset$. The problem is wide  open for $q\ge5$.

In a recent paper \cite{llp}, the second, the third, and the fourth authors managed to determine, as the first time,  the exact value of $r(C_4,B_n^{(2)})$ for infinitely many $n$. The  following upper bound was proved in \cite{llp}.
\begin{theorem}[Li, Lin, and Peng \cite{llp}]\label{llpub}
For all integers $q \ge4$ and $0 \le t\le q-1$, we have
\[
r(C_4,B_{(q-1)^2+(t-2)}^{(2)})\le q^2+t.
\]
\end{theorem}

We remark that the upper bound in Theorem \ref{llpub} improves that in \cite{FRS78}, see \eqref{frs}, by one for $q \ge4$ and $2 \le t\le q-2$. The improvement of the upper bound allows authors from \cite{llp} to prove the following result which gives the exact value of $r(C_4,B_n^{(2)})$ for infinitely many $n$.
\begin{theorem}[Li, Lin, and Peng \cite{llp}] \label{llplb}
For $q\geq 4$ being an even prime power, if $0 \le t \le q-1$ and $t\ne 1$, then
\[
  r(C_4,B_{(q-1)^2+(t-2)}^{(2)})= q^2+t.
\]
For $q\ge5$ being an odd prime power,  if $q\equiv3\pmod4$, $\tfrac{q+1}{2}\le t\le q-1$ and $t\neq \tfrac{q+3}{2}$; or $q\equiv1\pmod4$, $\tfrac{q-1}{2}\le t\le q-1$ and $t\neq \tfrac{q+1}{2}$, then
\[
  r(C_4,B_{(q-1)^2+(t-2)}^{(2)}) = q^2+t.
\]
\end{theorem}
In the same paper \cite{llp}, authors also established $r(C_4,B_n^{(2)}) \geq n+2\lfloor\sqrt{n}-6n^{0.2625}\rfloor$
for all large $n$. For more results on the Ramsey number of $C_m$ versus $B_n^{(k)}$ with $m \neq 4$, one can see, for instance \cite{ALPZ,FL,frs,HLLNP,LP21,ll,rs,shi}.

The goal of this paper is to extend Theorem \ref{llpub} and Theorem \ref{llplb} to all $k \geq 3$. First, we  prove the following upper bound for each $k \geq 3$, which is similar to Theorem \ref{llpub}. In the rest of this paper, for each integer $k \geq 3$  and a constant  $0<\varepsilon<1$, we define   $Q(k,\varepsilon)=(320k^4)^{k+1}/\varepsilon^{2k}$.

\begin{theorem}\label{mainub}
For each integer $k \geq 3$  and a constant  $0<\varepsilon<1$, if $q \geq Q(k,\varepsilon)$ and $0 \leq t \leq (1-\varepsilon)q$, then
\[
r(C_4,B^{(k)}_{q^2-kq+t+a_k}) \leq q^2+t,
\]
where $a_k=\binom{k}{2}-k$.
\end{theorem}

  Relying on constructions in \cite{llp}, we have the following result which gives  the exact value of $r(C_4,B_n^{(k)})$ for each $k \geq 3$ and infinitely many $n$.

\begin{theorem} \label{mainlb}
For each integer $k \geq 3$ and a constant $0 < \varepsilon <1$, let $a_k=\binom{k}{2}-k$.
For $q\geq Q(k,\varepsilon)$ being an even prime power, if $0 \le t \le (1-\varepsilon)q$ and $t\ne 1$, then
\[
 r(C_4,B^{(k)}_{q^2-kq+t+a_k})= q^2+t.
\]
For $q\geq  Q(k,\varepsilon)$ being an odd prime power,  if $q\equiv3\pmod4$, $\tfrac{q+1}{2}\le t\le (1-\varepsilon)q$ and $t\neq \tfrac{q+3}{2}$; or $q\equiv1\pmod4$, $\tfrac{q-1}{2}\le t\le (1-\varepsilon)q$ and $t\neq \tfrac{q+1}{2}$, then
\[
 r(C_4,B^{(k)}_{q^2-kq+t+a_k})= q^2+t.
\]
\end{theorem}
We do not attempt to optimize the constant $Q(k,\varepsilon)$.  In addition, we establish the following general upper and lower bounds for $r(C_4,B^{(k)}_n)$. Let $g_0(n)=n$ and $g_k(n)=g_{k-1}+\lfloor\sqrt{g_{k-1}-1}\rfloor+2$ for $k\ge 1$. Recall that we already know $r(C_4,B_{n}^{(k)}) \leq g_k(n)$ for $i \in \{1,2\}$.

\begin{theorem} \label{genub}
 For any positive integers $k \geq 3$ and sufficiently large $n$, we have
$$
 n+k\lfloor\sqrt n-6n^{0.2625}\rfloor-\frac 12k^2+\frac 32k \leq r(C_4,B_{n}^{(k)}) \le g_k(n).
$$
\end{theorem}

\noindent
{\bf Remark 1:} Note that  $g_k(n) \le n+k\lfloor\sqrt {n}\rfloor+\tfrac14(k^2+9k)$ by a rough computing. For $n=q^2-kq+t+a_k$ with $a_k=\binom{k}{2}-k$ and $0 \leq t \leq (1- \varepsilon)q$, one can see $n+k\lfloor\sqrt {n}\rfloor=q^2+t+\Theta(k^2)$.
Therefore, by Theorem \ref{mainlb}, the upper bound from Theorem \ref{genub} is tight up to a $\Theta(k^2)$ additive term.

\medskip
The organization of this paper is as follows. 
We will prove Theorem \ref{mainub} and Theorem \ref{mainlb} in Section \ref{main-t}. The proof of Theorem \ref{genub} will be presented in Section \ref{gen-lp}.

\medskip
\noindent
{\bf Notation:} For a graph $G=(V,E)$ with vertex set $V$ and edge set $E$, we use $\overline{G}$ to denote the complement of $G$. Given a subset $S \subset V(G)$, we write $G[S]$ as the subgraph of $G$ induced by $S$. 
For a vertex $v\in V$, let $N_G(v)$ be the neighborhood  of $v$ in $G$ and  $d_G(v)=|N_G(v)|$ be the degree of  $v$. We always use $\delta(G)$ and $\Delta(G)$ to denote the minimum degree and maximum degree in $G$ respectively.
  The subscription $G$ will be dropped if there is no confusion from the context.  We  omit the ceiling and floor if they do not affect the result. Let $[n]=\{1,2,\dots,n\}$ for a positive integer $n$.



\section{Proofs of Theorem \ref{mainub} and Theorem \ref{mainlb} }\label{main-t}
\subsection{Proof of Theorem \ref{mainub}}
To simplify the notation, let $n=q^2-kq+t+a_k$ and $b_k=a_k-\lceil \tfrac{k}{2}\rceil+2$, where $a_k=\binom{k}{2}-k.$  For $1 \leq i \leq k$, we define
\begin{equation*}
N_i=
\begin{cases}
  q^2-(k-i)q+t+b_k, &\textrm{ if } 1 \leq i \leq k-2; \\
  q^2-q+t+b_k+1,  & \textrm{ if } i=k-1; \\
   q^2+t,            & \textrm{ if } i=k.
  \end{cases}
\end{equation*}

We aim to show that for each integer $k \geq 3$  and a constant $0<\varepsilon<1$,  if $q \geq Q(k,\varepsilon)$ and $0 \leq t \leq (1-\varepsilon)q$, then
\[
r(C_4,B^{(k)}_{n}) \leq N_k.
\]


%

Let us recall the following classical result, which can be easily checked by using the double-counting method.
\begin{lemma}[K\H{o}v\'{a}ri, S\'{o}s, and Tur\'{a}n \cite{kst}]\label{C4-eg}
Let $G$ be a graph. If $C_4\nsubseteq G$, then $$\sum_{v\in V(G)} {d(v)\choose 2}\le{|V(G)| \choose 2}.$$

\end{lemma}

We start to prove the following lemma.
\begin{lemma}\label{lemma1}
For each $1 \leq i \leq k-1$ and  $q \geq Q(k,\varepsilon)$, we have $r(C_4,B_n^{(i)}) \leq N_i.$
\end{lemma}
\noindent
{\bf Proof.} The proof for $1 \leq i \leq k-2$ leverages on the induction on $i$. Note that $B_n^{(1)}$ is the star $K_{1,n}$. By Theorem \ref{pars}, we have that
\[
r(C_4,B_n^{(1)}) \leq n+\lfloor \sqrt{n-1} \rfloor+2=(q^2-kq+t+a_k)+(q-\lceil \tfrac{k}{2}\rceil)+2=N_1,
\]
here note that $ q-\lceil \tfrac{k}{2}\rceil \leq  \sqrt{n-1}  < q-\lceil \tfrac{k}{2}\rceil+1$ as $q$ is sufficiently large. For the induction step, let $G$ be a graph with $N_i$ vertices such that
$C_4 \not \subseteq G$ and $B_n^{(i)} \not \subseteq \overline{G}$. We next find a contradiction which implies that such a graph $G$ does not exist and the induction is complete.
Observe that $\Delta(\overline{G}) \leq N_{i-1}-1$. Otherwise, let $v$
be a vertex such that $|N_{\overline{G}}(v)| \geq N_{i-1}$. Let $S=N_{\overline{G}}(v)$.
 By the induction, either $G[S]$ contains a $C_4$ or $\overline{G}[S]$ contains a $B_n^{(i-1)}$. In the latter case, the vertex $v$ together with $B_n^{(i-1)}$ forms a $B_n^{(i)}$ in $\overline{G}$. There is a contradiction to the assumption on $G$ in each case, which implies that $\Delta(\overline{G}) \leq N_{i-1}-1$. Equivalently, $\delta(G) \geq N_i-N_{i-1}=q$.
 Recall $i \leq k-2$ and the assumption $t \leq (1-\varepsilon)q$ for some constant $\varepsilon$. Thus
 \[
 N_i= q^2-(k-i)q+t+b_k \leq q^2-2q+t+b_k \leq q^2-(1+\varepsilon)q+b_k < q^2-q,
 \]
provided $q \geq Q(k,\varepsilon)$.
 Now,
 \begin{equation} \label{eqgen}
 \sum_{v \in V(G)} \binom{d_G(v)}{2} \geq N_i \binom{q}{2}>\binom{N_i}{2}.
 \end{equation}
This contradicts Lemma \ref{C4-eg} and the induction is complete.

For $i=k-1$,  we can show $\delta(G) \geq N_{k-1}-N_{k-2}=q+1$ similarly.  Observe that inequality \eqref{eqgen} still holds, which also contradicts Lemma \ref{C4-eg}.  The proof of the lemma is complete. \hfill $\square$


\medskip

To prove Theorem \ref{mainub}, we need the following simple but novel refinement of Lemma \ref{C4-eg}.

\begin{lemma}\label{C4-eg1}
Let $G$ be a graph. If $C_4 \not \subseteq G$ and there are  $p$ pairs of vertices which are not
connected by a $2$-path (a path of length two),  then  $\sum_{v\in V(G)} {d(v)\choose 2} \leq {|V(G)| \choose 2}-p.$
\end{lemma}
\noindent
{\bf Proof.} Let $W=\{(\{x,y\}, z): x,y,z \in V(G) \textrm{ and } xz,yz \in E(G)\}$. Note that if $(\{x,y\}, z) \in W$, then $x$ and $y$ are connected by a 2-path.  In addition, for fixed  $x$ and $y$, there is at most one  $z$ such that $(\{x,y\}, z) \in W$ as $G$ is $C_4$-free.   Let $P \subset \binom{V(G)}{2}$ be the set of pairs of vertices which are not connected by a 2-path and $p=|P|$. Therefore, for $\{x,y\} \in P$, there is no $z$ such that $(\{x,y\}, z) \in W$. Utilizing the simple double-counting method, we get
\[
\sum_{v\in V(G)} {d(v)\choose 2}= |W| \leq {|V(G)| \choose 2}-p.
\]
The lemma is proved. \hfill $\square$

\medskip
We are ready to prove Theorem \ref{mainub}.

\medskip
\noindent
{\bf Proof of Theorem \ref{mainub}.}
 As we did in the proof of  Lemma \ref{lemma1}, let $G$ be a graph with $N_{k}$ vertices such that
$C_4 \not \subseteq G$ and $B_n^{(k)} \not \subseteq \overline{G}$. We next find a contradiction to this assumption, which will establish the upper bound.

 \vspace{0.2cm}
\noindent
{\bf Claim 1:}  $q-b_k-1 \leq \delta(G) \leq q$ and $\Delta(G) \leq q+b_k+2.$

\medskip
\noindent
{\bf Proof.} One can observe $\delta(G) \geq N_{k}-N_{k-1}=q-b_k-1$ similarly.
 Lemma \ref{C4-eg} yields that $\delta(G) \leq q$.  It is left to show $\Delta(G) \leq q+b_k+2.$ Otherwise, we pick a vertex $v$ such that $d(v):=d=\Delta(G)$. Let $N(v)=\{v_1,\dots,v_d\}$ and $A_i=N(v_i)\setminus(N(v)\cup\{v\})$. The assumption $C_4\nsubseteq G$ implies that edges in $N(v)$ must induce a matching in $G$ and $A_i\cap A_j=\emptyset$ for $1 \leq i\ne j \leq d$. It follows  that
 \[
 |A_i|\ge d(v_i)-2 \geq \delta(G)-2 \geq q-b_k-3.
 \]
Therefore, we have that
\[
N_{k}=|V(G)|\ge |\{v\}|+|N(v)|+\left|\bigcup_{i=1}^dA_i\right|\ge 1+d+d(q-b_k-3) \ge q^2+q-(b_k+3)^2,
\]
 here note that $d \geq q+b_k+3$ by the assumption.
However, since $t \leq (1-\varepsilon)q$, we know that
\[
N_{k}=q^2+t \leq q^2+(1-\varepsilon)q<q^2+q-(b_k+3)^2,
\]
where the last inequality holds as $q \geq Q(k,\varepsilon)$.
This is a contradiction and the claim is proved.  \hfill $\square$

 \vspace{0.2cm}
  Pick a vertex $v$ with degree $q-b_k-1 \leq d \leq q$ and define sets $A_i$ for $1 \leq i \leq d$ as above.  Notice that $A_i$ and $A_j$ are disjoint for $1\leq i \neq j \leq d$. Additionally, $q-b_k-3 \leq |A_i| \leq q+b_k+1$ since $q-b_k-1 \leq \delta(G) \leq \Delta(G) \leq q+b_k+2$.   It follows that $|\cup_{i=1}^d A_i  | \geq d(q-b_k-3) \geq (q-b_k-3)^2$  and
  $$\left|V(G)\setminus \left(\bigcup_{i=1}^d A_i\right) \right| \leq N_{k}-(q-b_k-3)^2=t+2(b_k+3)q+(b_k+2)^2  \leq (2b_k+7)q,$$
  where the last inequality holds as $q \geq Q(k,\varepsilon)$.
Let
\[
V_{\leq q}=\{v \in V(G):  d_G(v) \leq q\} \textrm{ and } V_{>q}=\{v \in V(G):  d_G(v) \geq q+1\}.
\]
As $G$ is $C_4$-free,  Lemma \ref{C4-eg} gives that
\begin{align*}
\binom{N_k}{2} &\geq \sum_{x \in V(G)} \binom{d_G(x)}{2}\\
                  &=\sum_{x \in V_{>q}} \binom{d_G(x)}{2}+\sum_{x \in V_{\leq q}} \binom{d_G(x)}{2}\\
                  & \geq \sum_{x \in V_{>q}} \binom{q+1}{2}+\sum_{x \in V_{\leq q}} \binom{q-b_k-1}{2}\\
                  &=|V_{>q}|\binom{q+1}{2}+|V_{\leq q}|\left(\binom{q+1}{2}-\binom{b_k+2}{2}-(b_k+2)(q-b_k-1) \right)\\
                  &\geq N_k\binom{q+1}{2}-(b_k+2)q|V_{\leq q}|.
\end{align*}
Solving the inequality above, we get that  $|V_{\leq q}|  \geq \eta|V(G)|$ where $\eta=\frac{\varepsilon}{2(b_k+2)}$.  As we already showed  $|V(G)\setminus \left(\cup_{i=1}^d A_i \right)|  \leq (2b_k+7)q$,  it follows that $| V_{\leq q} \cap \left(\cup_{i=1}^d A_i \right)  |  > \tfrac{\eta}{2}|V(G)|$ as $q \geq Q(k,\varepsilon)$.  Define
\[
J=\left\{1 \leq j \leq d: |A_j \cap V_{\leq q}| \geq \frac{\eta q}{8}\right\}.
\]
We assert that $|J| \geq \frac{\eta q}{8}$. Otherwise,
\begin{align*}
\left| V_{\leq q} \cap \left( \cup_{i=1}^d A_i \right)  \right|&=\sum_{i=1}^d |V_{\leq q} \cap A_i|=\sum_{i \in J}    |V_{\leq q} \cap A_i|+ \sum_{i \in [d] \setminus J}    |V_{\leq q} \cap A_i| \\
                                                 &\leq (q+b_k+1)|J|+\frac{\eta q}{8}(q-|J|)  \\
                                                 & \leq \frac{\eta q^2}{4}+\frac{\eta q}{8}(q+b_k+1)<\frac{\eta|V(G)|}{2},
\end{align*}
for  $q \geq Q(k,\varepsilon)$, which is a contradiction.

 For a $k$-set $K=\{x_1,\ldots,x_k\} \subset \cup_{i=1}^d A_i$,   we say $K$ is {\it admissible} if it satisfies the following properties:

\medskip
\noindent
$(1)$  $K$ is an independent set;

\noindent
$(2)$ $x_i \in V_{\leq q}$ for each $1 \leq i \leq k$;

\noindent
$(3)$  there is a $k$-set $\{i_1,\ldots,i_k\} \in \binom{[d]}{k}$ such that $|K \cap A_{i_j}|=1$ for each $1 \leq j \leq k$;

\noindent
$(4)$ there are no three vertices from $K$ with a common neighbor.

\medskip
Let $\cal K$ be the set of admissible $k$-sets. We have the following estimate for the size of $\cal K$.

\vspace{0.2cm}
\noindent
{\bf Claim 2:} $|{\cal K}| \geq c_{k,\varepsilon}q^{2k}$ for some constant $c_{k,\varepsilon} >0$.

\vspace{0.2cm}
\noindent
{\bf Proof.} Recall the definition of $J$. As  $|J| \geq \frac{\eta q}{8}$, we may assume $|A_i \cap V_{\leq q}| \geq \frac{\eta q}{8}$ for each $1 \leq i \leq \frac{\eta q}{8}$. In addition, let $B_i \subseteq A_i \cap V_{\leq q}$ with $|B_i|= \frac{\eta q}{8}$ for each $1 \leq i\le \frac{\eta q}{8}$. We next construct admissible $k$-sets greedily. We first select a $k$-set of indices $\{i_1,\ldots,i_k\} \in \binom{[\eta q/8]}{k}$. For a given index set,  we pick $x_1 \in B_{i_1}$ and $x_2 \in B_{i_2}\backslash N_G(x_1)$ arbitrarily and  select vertices $x_j$ from $B_{i_j}$ one by one recursively for $3 \leq j \leq k$.
Observe that if we choose $k$ vertices in above manner, then the resulting $k$-set satisfies Properties (2) and  (3).
Clearly, $|B_{i_2}\backslash N_G(x_1)|\ge|B_{i_2}|-1$ as $G$ is $C_4$-free. Assume that we already choose vertices $x_1,\ldots,x_j$ such that $\{x_1,\ldots,x_j\}$ satisfies Properties (1) and  (4). To make the notation simple, let $S=\{x_1,\ldots,x_j\}$.
We next show how to pick $x_{j+1} \in B_{i_{j+1}}$ so that the set $S \cup \{x_{j+1}\}$ maintains Properties $(1)$ and $(4)$. Note that $A_i=N(v_i)\setminus(N(v)\cup\{v\})$ for each $1 \leq i \leq d$. The assumption $C_4 \not \subseteq G$ implies that each vertex $y \in V(G)\setminus \{v_i\}$ has at most one neighbor in $A_i$ and of course in $B_i$.
Let $$B_{i_{j+1}}'=B_{i_{j+1}} \setminus \left(\cup_{\ell=1}^j N_G(x_{i_\ell})\right).$$ It follows that $|B_{i_{j+1}}'| \geq \frac{\eta q}{8}-j \geq \frac{\eta q}{8}-k$.  Note that  the set $S \cup \{x_{j+1}\}$ satisfies Property $(1)$ for each $x_{j+1} \in B_{i_{j+1}}'$.

For each  2-set $\{y,z\} \in \binom{S}{2}$,  let $w_{yz}$ be the possible common neighbor of $y$ and $z$ in $G$. Note that $w_{yz} \neq v_{i_{j+1}}$ since $y \not \in A_{i_{j+1}}$ and $z \not \in A_{i_{j+1}}$.
As $G$ is $C_4$-free, $w_{yz}$ has at most one neighbor in $B_{i_{j+1}}'$. Let
\[
B_{i_{j+1}}''=B_{i_{j+1}}'\setminus \left( \cup_{\{y,z\} \in \binom{S}{2}} N_G(w_{yz})\right).
\]
Thus
\[
|B_{i_{j+1}}''| \geq |B_{i_{j+1}}'|-\binom{j}{2} \geq \frac{\eta q}{8}-k-\binom{k}{2} \geq \frac{\eta q}{10}
\]
 since $q \geq Q(k,\varepsilon)$.  Recall that no three vertices from $S$ have a common neighbor.
  We assert that $S \cup \{x_{i_{j+1}}\}$ satisfies Property (4) for each  $x_{i_{j+1}} \in B_{i_{j+1}}''$.
 Otherwise, suppose there are three vertices from  $S \cup \{x_{i_{j+1}}\}$ with a common neighbor. As $S$ satisfies Property (4) by the assumption, we may assume  $y,z,$ and  $x_{i_{j+1}}$ having a common neighbor $w_{yz}$, where $y,z \in S$.
As we removed vertices from $N_G(w_{yz})$ in the course of defining $B_{i_{j+1}}''$, the vertex $x_{i_{j+1}}$ is not a neighbor of $w_{yz}$. This is  a contradiction and the assertion follows.
  Therefore,
  the set $S \cup \{x_{j+1}\}$ satisfies Properties (1) and (4)
 for each $x_{j+1} \in B_{i_{j+1}}''$ by noting $B_{i_{j+1}}' \subseteq B_{i_{j+1}}''$.
  Notice that once we finish the selection of the vertex $x_{k}$ in the above way. The set $\{x_1,\ldots,x_k\}$ is an admissible set. Apparently, there are  $\binom{\eta q/8}{k}$ choices of the index set
and there are at least $(\eta q/10)^k$ choices of vertices $x_1,\ldots,x_k$ for a given index set.  Therefore,
\[
|{\cal K}| \geq \left(\frac{\eta q}{10}\right)^k \binom{\eta q/8}{k}=c_{k,\varepsilon}q^{2k}.
\]
The proof is complete. \hfill $\square$

\vspace{0.2cm}
For a 2-set of vertices $P=\{u,v\}$, we say $P$ is a {\it good pair} if $N_G(u) \cap N_G(v)=\emptyset$ and $\{u,v\} \subset V_{\leq q}$.

\vspace{0.2cm}
\noindent
{\bf Claim 3:} \ For each admissible $k$-set $K=\{x_1,\ldots,x_k\}$, there are two vertices $x_i, x_j \in K$  such that $\{x_i,x_j\}$ is a good pair.

\vspace{0.2cm}
\noindent
{\bf Proof.}  As Property $(2)$, we need only to show that there are two vertices $x_i, x_j \in K$ with $N_G(x_i) \cap N_G(x_j)=\emptyset$. If such kind pair of vertices does not exist, then $|N_G(x_i) \cap N_G(x_j)|=1$ for any $1 \leq i < j \leq k$ because $G$ does not contain $C_4$. By the inclusion-exclusion formula and Property $(4)$, we get
\[
 \left | \cup_{i=1}^{k} N_G(x_i) \right|=\sum_{i=1}^k |N_G(x_i)|-\binom{k}{2} \leq kq-\binom{k}{2},
\]
where the last inequality follows from Property $(2)$. Let $Y=V(G) \setminus  ( \cup_{i=1}^k N_G[x_i])$, here $N_G[x_i]=N_G(x_i) \cup \{x_i\}$. Then  $K$ and $Y$ are completely adjacent in  $\overline{G}$.  Notice that $|Y| \geq N_k-k-kq+\binom{k}{2}=n$ and  $K$ is a clique in $\overline{G}$ as Property $(1)$. It follows that $B_n^{(k)} \subseteq \overline{G}$, which is a contradiction to the definition of $G$. The claim is proved. \hfill $\square$

\vspace{0.2cm}
Let $\cal P$ be the set of good pairs which is contained in some admissible set from ${\cal K}$.

\vspace{0.2cm}
\noindent
{\bf Claim 4:}  $|{\cal P}| \geq d_{k,\varepsilon}q^{4}$ for some constant $d_{k,\varepsilon} >0$.

\vspace{0.2cm}
\noindent
{\bf Proof.} Let $\{x_{i_1},x_{i_2}\}$ be a good pair with $x_{i_1} \in A_{i_1}$ and $x_{i_2} \in A_{i_2}$.
For a fixed index set $\{i_3,\ldots,i_k\}$, there are at most $\prod_{j=3}^k |A_{i_j}|$ choices for $\{x_{i_3},\ldots,x_{i_k}\}$ such that
 $$\{x_{i_1},x_{i_2},x_{i_3},\ldots,x_{i_k}\} \in {\cal K}.$$
There are at most $\binom{d-2}{k-2}$ choices for the index set $\{i_3,\ldots,i_k\}$. Recall $|A_{i_j}| \leq q+b_k+1$ and $d \leq q$. Therefore, each good pair is contained in at most $(q+b_k+1)^{k-2}\binom{q-2}{k-2}$ admissible $k$-sets from $\cal K$. As we already showed that  $|{\cal K}| \geq c_{k,\varepsilon}q^{2k}$, see Claim 2, it follows that
\[
|{\cal P}| \geq \frac{c_{k,\varepsilon}q^{2k}}{(q+b_k+1)^{k-2}\binom{q-2}{k-2}}>d_{k,\varepsilon}q^4,
\]
for some constant $d_{k,\varepsilon}>0$, which completes the proof. \hfill $\square$

\vspace{0.2cm}
Note that vertices from a good pair are not connected by a 2-path. Applying Lemma \ref{C4-eg1} with $p=d_{k,\varepsilon}q^4$, we get
\[
\frac{q^4}{2}-s_kq^3<N_k \binom{q-b_k-1}{2} \leq \sum_{v \in V(G)} \binom{d_G(v)}{2} \leq \binom{N_k}{2}-d_{k,\varepsilon}q^4 \leq (0.5-d_{k,\varepsilon})q^4+t_kq^3,
\]
where $s_k$ and $t_k$ are constants depending on $k$. As $q \geq Q(k,\varepsilon)$ is large enough, we obtain a contradiction which finishes the proof.
 \hfill $\square$

\subsection{Proof of Theorem \ref{mainlb}}
To prove Theorem \ref{mainlb}, it remains to show the lower bounds since the upper bounds have been established in Theorem \ref{mainub}. Actually, we will prove a slightly stronger result.

\vspace{0.2cm}
\noindent
{\bf Proof of Theorem \ref{mainlb}.}  For $q$ being an even prime power and $0 \le t \le q$ with $t\ne 1$, we next show
\[
 r(C_4,B^{(k)}_{q^2-kq+t+a_k}) \geq  q^2+t,
\]
where $a_k=\binom{k}{2}-k$.
For each $0 \le t \le q$ with $t\ne 1$,  authors of \cite{llp} constructed a $C_4$-free graph $H_q^t$ with $q^2+t-1$ vertices which is an induced subgraph of the well-known Erd\H{o}s-R\'enyi orthogonal polarity graph \cite{Br66,ER66}.  In addition, $\delta(H_q^t) \geq q$. Note that each pair of vertices has at most one common neighbor. It is sufficient to show the complement of $H_q^t$ does not contain $B^{(k)}_{q^2-kq+t+a_k}$ as a subgraph.
For a fixed $k$-subset $U$,  the inclusion-exclusion formula yields that
\begin{align*}
\left|\cup_{i=1}^k N(u_i) \right| & \geq \sum_{u \in U} N(u) -\sum_{\{u_i,u_j\} \in \binom{U}{2}} |N(u_i) \cap N(u_j)| \\
                                  & \geq kq-\binom{k}{2}.
\end{align*}
Therefore, each set of $k$ vertices in the complement of $H_q^t$ has at most $q^2+t-1-(kq-\binom{k}{2}+k)=q^2-kq+t+a_k-1$ common neighbors and the lower bound in this case follows.

For $q$ being an odd prime power,  if either $q\equiv3\pmod4$, $\tfrac{q+1}{2}\le t\le q$ and $t\neq \tfrac{q+3}{2}$, or $q\equiv1\pmod4$, $\tfrac{q-1}{2}\le t\le q$ and $t\neq \tfrac{q+1}{2}$, authors of \cite{llp} constructed a $C_4$-free graph $G_q^t$ with $q^2+t-1$ vertices and $\delta(G_q^t) \geq q$. One can prove that the complement of $G_q^t$ does not contain $B^{(k)}_{q^2-kq+t+a_k}$ as a subgraph similarly. The proof of Theorem \ref{mainlb} is complete. \hfill $\square$

\section{Proof of Theorem \ref{genub}}\label{gen-lp}

In this section, we give the proof of Theorem \ref{genub}. To prove the upper bound, we will need the following lemma by  Erd\H{o}s, R\'enyi and S\'os \cite{ER66}.
The {\em Friendship graph} or a {\em
$k$-fan} is the graph $F_k=K_1+kK_2$ which consists of $k$ triangles with a vertex in common.
\begin{lemma} \label{friend}
Let $G$ be a graph with $N$ vertices such that any pair of
vertices is joined by exactly one path of length two in $G$. Then
$N=2k+1$ and $G=F_k$.
\end{lemma}

\medskip
\noindent
{\bf Proof of the upper bound.} For simplicity, let $g_k=g_k(n)$ for $k\ge 1$. We apply induction on $k$  to show $r(C_4,B_n^{(k)}) \leq g_k$ for $k \geq 1$. The base case where $k=1$ was done by Parsons \cite{TP75}.
Now assume that $k\ge 2$ and the assertion holds for $k-1$, i.e., $r(C_4,B_{n}^{(k-1)}) \le g_{k-1}$.  We aim to show $r(C_4,B_{n}^{(k)}) \le g_{k}$.

Suppose that $G$ is a graph on $g_k$ vertices such that $G$ is $C_4$-free and $\overline{G}$ is $B_{n}^{(k)}$-free. It is clear that the maximum degree $\Delta(\overline{G})\le r(C_4,B_{n}^{(k-1)})-1\le g_{k-1}-1$, so we have
$
\delta(G)\ge g_k-g_{k-1}.
$
Since $G$ is $C_4$-free, Lemma \ref{C4-eg} yields that
\begin{equation}\label{parsons-1}
\sum_{w\in V(G)} {d(w)\choose 2} \le {g_k\choose 2}.
\end{equation}
As $\delta(G)\ge g_k-g_{k-1}$, it follows that
\begin{equation}\label{parsons-2}
g_k^2 -2(g_{k-1}+1)g_{k} +g_{k-1}^2+g_{k-1}+1 \le 0.
\end{equation}
The above inequality becomes an equality  if and only if $G$ is
$(g_k-g_{k-1})$-regular  and each pair of vertices have exactly
one common neighbor.

If the inequality (\ref{parsons-2})
is an equality, then  $G=F_n$ by Lemma \ref{friend}. However,  $G$ must be
regular, which is a contradiction as $n\ge 3$.  This is not the case and hence
the inequality (\ref{parsons-2}) is strict, which implies that
\[
g_k^2 -2(g_{k-1}+1)g_{k} +g_{k-1}^2+g_{k-1}+2 \le 0,
\]
here note that $g_k$ is an integer.
Solving the inequality above, we get
\[
g_k \leq g_{k-1}+ \sqrt{g_{k-1}-1} +1.
\]
Therefore, $g_k \leq g_{k-1}+ \lfloor \sqrt{g_{k-1}-1} \rfloor +1$ as $g_k$ is an integer and this is a contradiction to the definition of $g_k$.
 This completes the induction and hence the proof of the upper bound.
\hfill $\Box$


\medskip
\noindent
{\bf Proof of the lower bound.} Let $\alpha=0.525$ and $p$ be the smallest prime with $p\ge n^{1/2}+\frac12$.
  A result in \cite{bhp}  asserts that for sufficiently large $k$ there is a prime in the interval $(k,k+k^{0.525}]$.
 Thus  we get $p\le n^{1/2}+n^{\alpha/2}+1$. For each prime $p$, it is known, see \cite{ER66,Br66}, that there exists a $C_4$-free graph $ER_p$ of order $N = p^2 +p+ 1$  in which the degree of each vertex is $p$ or $p + 1$.
Set $m = \lfloor n^{1/2}-6n^{\alpha/2}\rfloor$ and $d = N-(n+mk-\frac 12k^2+\frac 32k-1)$.
Let us randomly delete $d$ vertices from $ER_p$ to obtain a random subgraph $G$. We say a vertex $v$ in $ER_p$ is ``bad'' if it is not deleted and has degree less than $m$ in $G$.
Let $B_v$ denote the event that $v$ is bad. Note that a vertex $v$  of degree $p$ is bad if and only if $v$ is not deleted and $t$  neighbors of $v$ are removed for  some $t> p-m$. Therefore,
\[
\Pr(B_v)=\sum_{t > p-m}\frac{{p\choose t}{{N-p-1}\choose {d-t}}}{{N\choose d}}
<\sum_{t > p-m}\left(\frac{epd}{t(N-d)}\right)^t<N\left(\frac{2e}{6}+o(1)\right)^{6n^{\alpha/2}}
\]
by noting $p\le n^{1/2}+n^{\alpha/2}+1$, $d= N-(n+mk-\frac 12k^2+\frac 32k-1)\le (2+o(1))n^{(1+\alpha)/2}$, $t> p-m > 6n^{\alpha/2}$ and $N-d=n+mk-\frac 12k^2+\frac 32k-1>n$ for any fixed $k$ and large $n$.  Clearly, a similar argument works for a vertex $v$ with degree $p+1$. Let $X$ be the random variable counting the number of bad vertices. By the linearity of expectation,  E$(X)=\sum_{v \in V(ER_p)}\Pr(B_v)$. As the upper bound for $\Pr(B_v)$, we get  E$(X) <1$ for   large $n$. Thus there is a choice of $d$ vertices such that the corresponding subgraph $G$ contains no bad vertices.  Equivalently, there exists a $C_4$-free graph $G$ on $n+k\lfloor\sqrt n-6n^{0.2625}\rfloor-\frac 12k^2+\frac 32k-1$ vertices with minimum degree at least $m = \lfloor n^{1/2}-6n^{\alpha/2}\rfloor$.

 We next show $\overline{G}$ is $B_n^{(k)}$-free.  Let $v_1,\ldots,v_k \in V(G)$.  Observe  that $|N_G(v_i)\cap N_G(v_j)|\le 1$ for $1\le i\ne j \le k$ since $C_4\nsubseteq G$. By the inclusion-exclusion formula, we obtain that
\[
\displaystyle \left|\bigcup\limits_{i=1}^k N(v_i)\right|\ge  \sum_{i=1}^k |N_G(v_i)|-\sum_{1\leq i <j \leq k}|N_G(v_i) \cap N_G(v_j)| \geq km-\binom{k}{2}=\frac{1}{2}k(2m-k+1),
\]
the last inequality follows from $\delta(G) \geq m$.
 Therefore, the number of common neighbors of any $k$ vertices in $\overline{G}$ is at most
 \begin{align*}
 |V(G)|-k-\left|\bigcup\limits_{i=1}^k N(v_i)\right|\le \left(\displaystyle n+mk-\frac 12k^2+\frac 32k-1\right)-k-\frac 12k(2m-k+1)=n-1.
 \end{align*}
 Therefore, $G$ is $C_4$-free and $\overline{G}$  is $B_n^{(k)}$-free.  It follows that $\displaystyle r(C_4,B_n^{(k)})> n+mk-\frac 12k^2+\frac 32k-1$ and the lower bound is proved.
\hfill$\Box$

%
%

\end{spacing}

\end{document}